\documentclass{amsart}
\usepackage{amsmath,amsthm,amssymb,amsfonts}

\theoremstyle{plain}
\newtheorem{theo+}           {Theorem}      [section]
\newtheorem{prop+}  [theo+]  {Proposition}
\newtheorem{lemm+}  [theo+]  {Lemma}
\newtheorem{cor+}  [theo+]  {Corollary}
\newenvironment{theorem}{\begin{theo+}}{\end{theo+}}

\newenvironment{corollary}{\begin{cor+}}{\end{cor+}}

 \newcommand{\De}{\Delta}

\begin{document}
\baselineskip 18pt
\larger[2]
\title[New transformations for elliptic hypergeometric series]
{New transformations for elliptic\\
  hypergeometric series  on the root system $A_n$}
\author{Hjalmar Rosengren}
\address{Department of Mathematics\\ Chalmers University of Technology and 
G\"oteborg University\\SE-412~96 G\"oteborg, Sweden}
\email{hjalmar@math.chalmers.se}
\urladdr{http://www.math.chalmers.se/{\textasciitilde}hjalmar}

\begin{abstract}
Recently, Kajihara gave a Bailey-type transformation relating
basic hypergeometric series on the root system $A_n$, with different
dimensions $n$. We give, with a new, elementary, proof, an
 elliptic analogue of this transformation. We also obtain  further
Bailey-type transformations as consequences of our result, some of
which are new also in the case of basic and classical hypergeometric
series. 
\end{abstract}
\maketitle   

\section{Introduction}

Elliptic hypergeometric series form an extension of classical and basic
(or $q$-) hyper\-geometric series, which was introduced by 
Frenkel and Turaev in 1997 \cite{ft}. It was found that Jackson's 
${}_8W_7$ summation and Bailey's ${}_{10}W_9$ transformation 
  admit one-parameter extensions, roughly speaking obtained
by replacing ``$1-x$'' by the theta function
$\prod_{j=0}^\infty(1-p^jx)(1-p^{j+1}/x)$. For elliptic
hypergeometric series, the so called balanced  and well-poised
conditions on the series appearing in these identities reflect
invariance properties under the modular group \cite{ft,s}.

In the last few years, multivariable elliptic hypergeometric series, in
particular series associated to classical root systems, has received
much attention \cite{ds1,ds2,rw,r2,rs,s,sp,w}. In
the present paper we  build upon the work in \cite{r2} to obtain
some new transformation formulas for elliptic hypergeometric series on
the root system $A_n$. 

In Theorem \ref{t} we  give an elliptic analogue of a multivariable Bailey
transformation recently
discovered by Kajihara \cite{ka}. In contrast to 
most known transformations,
Kajihara's identity relates sums of \emph{different} dimension; 
see\ \cite{gkr,kr,r,rc} for further results with this
property.  (We mention that, in view of the
analogy  between
 hypergeometric series and hypergeometric integrals, there may exist related
transformations between integrals of different dimension. The only such
 result we are aware of is in the recent paper \cite{tv}; this seems
  not directly related to series of the type studied here, 
but rather to discrete Selberg integrals \cite{a}.)

In Section 4 we obtain further  Bailey-type 
transformations, between series of \emph{the same}
dimension, by iterating  Theorem \ref{t}. Most of these are new
also in the case of basic and classical hypergeometric series.

\section{Notation}

Elliptic hypergeometric series may be built from the theta function
\begin{equation}\label{th}
\theta(x)=\prod_{j=0}^\infty(1-p^jx)(1-p^{j+1}/x),\qquad |p|<1.\end{equation}
We will often use that
$\theta(1/x)=-\theta(x)/x.$

We denote elliptic Pochhammer symbols by
$$(a)_k=\theta(a)\,\theta(aq)\dotsm\theta(aq^{k-1}).$$
The constants $p$ and $q$ are fixed throughout the paper and will be suppressed
from the notation. The elementary identities
$$(a)_{n+k}=(a)_n(aq^n)_k,$$
$$(a)_{n-k}=(-1)^kq^{\binom k2}(q^{1-n}/a)^k
\frac{(a)_n}{(q^{1-n}/a)_k},$$
$$(a)_n=(-1)^nq^{\binom n2}a^n(q^{1-n}/a)_n $$
 will be used repeatedly and without comment.
Occasionally, we  use the shorthand notation
$$(a_1,\dots,a_n)_k=(a_1)_k\dotsm(a_n)_k.$$

We write
$$\De(z)=\prod_{1\leq j<k\leq n} z_j\,\theta(z_k/z_j).$$
This may be viewed as an elliptic analogue of the Weyl denominator for
the root system $A_n$. 
Elliptic  hypergeometric series on  $A_n$ are
characterized by the factor
$$\frac{\De(zq^y)}{\De(z)}=\prod_{1\leq j<k\leq n}q^{y_j}
\frac{\theta(z_kq^{y_k}/z_jq^{y_j})}{\theta(z_k/z_j)}, $$
where the $z_k$ are free parameters and the $y_k$ summation indices.

Note that when $p=0$ we have
$$(a)_k=(1-a)(1-aq)\dotsm(1-aq^{k-1}), $$
the standard building blocks of basic hypergeometric series. 
Moreover, in this case, 
$$\frac{\De(zq^y)}{\De(z)}=\prod_{1\leq j<k\leq n}
\frac{z_kq^{y_k}-z_jq^{y_j}}{z_k-z_j}. $$
Rescaling and letting $q\rightarrow 1$ one recovers the classical
Pochhammer symbols
$$a(a+1)\dotsm(a+k-1)$$
 and double product
$$\frac{\Delta(z+y)}{\Delta(z)}=\prod_{1\leq j<k\leq n}
\frac{z_k+y_k-z_j-y_j}{z_k-z_j}, $$
characterizing classical hypergeometric series on $A_n$ \cite{hbl}.

For later reference we give some facts about one-variable elliptic
hypergeometric series.
Let $E$ be the function
\begin{equation*}\begin{split}&
E(a;q^{-N},b,c,d,e,f,g)\\
&\quad=\sum_{k=0}^N\frac{\theta(aq^{2k})}{\theta(a)}
\frac{(a,q^{-N},b,c,d,e,f,g)_k}{(q,aq^{N+1},aq/b,aq/c,aq/d,aq/e,aq/f,aq/g)_k}
\,q^k.\end{split}\end{equation*}
This is a ${}_{10}W_9$ sum when $p=0$. Frenkel and Turaev \cite{ft} proved the
transformation formula
\begin{equation}\label{fb}\begin{split}
&E(a;q^{-N},b,c,d,e,f,g)\\
&\quad=\frac{(aq,aq/ef,\lambda q/e,\lambda
  q/f)_N}{(aq/e,aq/f,\lambda q,\lambda q/ef)_N}\,
E(\lambda;q^{-N},\lambda b/a,\lambda c/a,\lambda d/a,e,f,g),
\end{split}\end{equation}
where $bcdefg=a^3q^{N+2}$ and $\lambda=qa^2/bcd$.
When $p=0$, this is the famous Bailey transformation \cite{b,gr}.       
Iterating \eqref{fb} one obtains 
\begin{equation}\label{sb}\begin{split}&
E(a;q^{-N},b,c,d,e,f,g)=g^N\frac{(aq/cg,aq/dg,aq/eg,aq/fg,aq,b)_N}
{(aq/c,aq/d,aq/e,aq/f,aq/g,b/g)_N}\\
&\quad\times E(gq^{-N}/b;q^{-N},gq^{-N}/a,
aq/bc,aq/bd,aq/be,aq/bf,g),\end{split}\end{equation}
still assuming  $bcdefg=a^3q^{N+2}$. When $p=0$, this is
\cite[Exercise 2.19]{gr}.

\section{An elliptic Kajihara transformation}

The following identity is our main result.

\begin{theorem}\label{t}
Assuming 
\begin{equation}\label{bc}
w_1\dotsm w_m=z_1\dotsm z_na_1\dotsm a_{m+n},\end{equation}
 the following identity holds:
\begin{equation}\label{te}\begin{split}&
\sum_{\substack{y_1,\dots,y_n\geq 0\\y_1+\dots+y_n=N}}
\frac{\De(zq^y)}{\De(z)}\prod_{k=1}^n\frac{\prod_{j=1}^{m+n}(a_jz_k)_{y_k}}
{\prod_{j=1}^m(w_jz_k)_{y_k}\prod_{j=1}^n(qz_k/z_j)_{y_k}}\\
&\qquad=\sum_{\substack{y_1,\dots,y_m\geq 0\\y_1+\dots+y_m=N}}
\frac{\De(wq^y)}{\De(w)}
\prod_{k=1}^m\frac{\prod_{j=1}^{m+n}(w_k/a_j)_{y_k}}
{\prod_{j=1}^n(w_kz_j)_{y_k}\prod_{j=1}^m(qw_k/w_j)_{y_k}}.
\end{split}\end{equation}
\end{theorem}

The case $m=n=2$ is easily seen to be equivalent to
\eqref{sb}, so Theorem \ref{t} is a multivariable generalization
of this transformation. On the other hand, if
 $m=1$ but $n$ is general we have the summation formula
\begin{equation}\label{j}
\sum_{\substack{y_1,\dots,y_n\geq 0\\y_1+\dots+y_n=N}}
\frac{\De(zq^y)}{\De(z)}\prod_{k=1}^n\frac{\prod_{j=1}^{n+1}(a_jz_k)_{y_k}}
{(wz_k)_{y_k}\prod_{j=1}^n(qz_k/z_j)_{y_k}}=\frac{\prod_{j=1}^{n+1}(w/a_j)_{N}}
{\prod_{j=1}^n(wz_j)_{N}(q)_{N}},
\end{equation}
where  $w=z_1\dotsm z_na_1\dotsm a_{n+1}$.
This is 
 \cite[Theorem 5.1]{r2}, which is an elliptic analogue of 
Milne's  $A_n$ Jackson summation  \cite{m}.
 See also  \cite{sp}, where
 it was shown that \eqref{j} follows from a
certain conjectured multiple integral evaluation.

The non-elliptic case, $p=0$,
 of Theorem \ref{t} is equivalent to Proposition 6.2 of
\cite{ka}, where it was derived using
 Macdonald polynomials. A different proof,
based on Gustafson's $A_n$ ${}_6\psi_6$ summation \cite{gu},
was given in \cite{r}. Neither of these proofs is likely to admit a
straight-forward elliptic generalization. Here we will use
 a  simple inductive argument, which works also in the elliptic case.

\begin{proof}

We  prove Theorem \ref{t} by induction on $n$. As a starting point we
need the case $n=1$, or equivalently $m=1$, that is, the identity
\eqref{j}. 

Assume that Theorem \ref{t} holds for fixed $n$ but general $m$,
and denote the left-hand side of \eqref{te} by $S_{nm}(z;w;a)$.
Consider the sum  $S_{n+1,m}(z;w;a)$. We replace the index set
$(y_1,\dots,y_{n+1})$ by $(y_1,\dots,y_{n},s)$ and rewrite part of the
summand as
\begin{equation*}\begin{split}&
\frac{\De( zq^{(y,s)})}{\De(
  z)}\prod_{j=1}^{n+1}\frac{1}{(qz_{n+1}/z_j)_s}
=\frac{1}{(q)_s}\frac{\De(\tilde zq^y)}{\De(\tilde z)}\prod_{j=1}^n
\frac{q^{-y_j}\theta(z_{n+1}q^s/z_jq^{y_j})}{\theta(z_{n+1}/z_j)
(qz_{n+1}/z_j)_s }\\
&\quad=\frac{1}{(q)_s}\frac{\De(\tilde zq^y)}{\De(\tilde z)}\prod_{j=1}^n
\frac{\theta(z_jq^{y_j-s}/z_{n+1})}{(z_{n+1}/z_j)_s 
\theta(z_jq^{-s}/z_{n+1})}\\
&\quad=\frac{1}{(q)_s}\frac{\De(\tilde zq^y)}{\De(\tilde z)}\prod_{j=1}^n
\frac{(z_jq^{1-s}/z_{n+1})_{y_j}}{(z_{n+1}/z_j)_s 
(z_jq^{-s}/z_{n+1})_{y_j}},
\end{split} \end{equation*}
where   $\tilde z=(z_1,\dots,z_n)$. This gives
\begin{multline*}S_{n+1,m}(z;w;a)
=\sum_{s=0}^N\frac{\prod_{j=1}^{m+n+1}(a_jz_{n+1})_s}
{(q)_s\prod_{j=1}^m(w_jz_{n+1})_s\prod_{j=1}^n(z_{n+1}/z_j)_s}\\
\times\!\sum_{y_1+\dots+y_n=N-s}\!
\frac{\De(\tilde zq^y)}{\De(\tilde z)}
\prod_{k=1}^n\frac{(q^{1-s}z_k/z_{n+1})_{y_k}\prod_{j=1}^{m+n+1}(a_jz_k)_{y_k}}
{(q^{-s}z_k/z_{n+1})_{y_k}(qz_k/z_{n+1})_{y_k}
\prod_{j=1}^m(w_jz_k)_{y_k}\prod_{j=1}^n(qz_k/z_j)_{y_k}}.
 \end{multline*}

We  observe that the inner sum is of the form
$S_{n,m+2}(\tilde z;\tilde w;\tilde a)$, with
\begin{equation}\label{aw}\begin{split}
\tilde w&=(w_1,\dots,w_m,q/z_{n+1},q^{-s}/z_{n+1}),\\
\tilde a&=(a_1,\dots,a_{m+n+1},q^{1-s}/z_{n+1}).\end{split}\end{equation}
We apply the induction hypothesis to this sum,
 writing $y_{m+1}=t$,
$y_{m+2}=u$, and then change the order of
summation according to
\begin{equation}\label{co}\sum_{s=0}^N\,\sum_{y_1+\dots+y_m+t+u=N-s}(\dotsm)=
\sum_{t=0}^N\,\sum_{y_1+\dots+y_m\leq
  N-t}\,\sum_{s+u=N-t-|y|}(\dotsm).\end{equation} 

We consider first the inner sum, collecting all factors involving $s$
and $u$:
\begin{equation}\label{su}\begin{split}&
\frac{\prod_{j=1}^{m+n+1}(a_jz_{n+1})_s}
{(q)_s\prod_{j=1}^m(w_jz_{n+1})_s\prod_{j=1}^n(z_{n+1}/z_j)_s}
\prod_{j=1}^m\left(q^{y_j}\frac{\theta(\tilde w_{m+2}q^u/w_jq^{y_j})}
{\theta(\tilde w_{m+2}/w_j)}\frac{(w_j/\tilde
    a_{m+n+2})_{y_j}}{(qw_j/\tilde w_{m+2})_{y_j}}\right)\\
&\quad\times q^t
\frac{\theta(\tilde w_{m+2}q^u/\tilde w_{m+1}q^{t})}
{\theta(\tilde w_{m+2}/\tilde w_{m+1})}
\frac{(\tilde w_{m+1}/\tilde
    a_{n+m+2})_{t}}{(q\tilde w_{m+1}/\tilde w_{m+2})_{t}}
\frac{\prod_{j=1}^{m+n+2}(\tilde w_{m+2}/\tilde
    a_{j})_{u}}{\prod_{j=1}^n(\tilde
    w_{m+2}z_j)_u\prod_{j=1}^{m+2}(q\tilde w_{m+2}/\tilde w_j)_u}\\
&=q^t\frac{\theta(q^{u-s-t-1})}{\theta(q^{-s-1})}\frac{(q^s)_t(q^{-1})_u}{(q)_s
(q^{s+2})_t(q^{-s})_u(q)_u}
\frac{\prod_{j=1}^{m+n+1}(a_jz_{n+1})_s(q^{-s}/z_{n+1}a_j)_u}
{\prod_{j=1}^n(z_{n+1}/z_j)_s(q^{-s}z_j/z_{n+1})_u}\\
&\quad\times\prod_{j=1}^m\left(q^{y_j}\frac{\theta(q^{u-s-y_j}/z_{n+1}w_j)}
{\theta(q^{-s}/z_{n+1}w_j)}\frac{(q^{s-1}w_jz_{n+1})_{y_j}}
{(w_jz_{n+1})_s(q^{1+s}w_jz_{n+1})_{y_j}(q^{1-s}/z_{n+1}w_j)_u}\right).
\end{split}\end{equation}

Note that, because of the factor
 $(q^{-1})_u$, \eqref{su} vanishes unless
$u\in\{0,1\}$. (If $s=0$ and $u=1$, the factor $1/(q^{-s})_u$
gives an apparent singularity, but this is removed by 
$(q^s)_t$ if $t>0$ and by $\theta(q^{u-s-t-1})$ if $t=0$.) This leads
 to considerable simplification. For instance,
\begin{equation*}\begin{split}& 
(a_jz_{n+1})_s(q^{-s}/z_{n+1}a_j)_u=(-1)^uq^{\binom
  u2}(q^{-s}/z_{n+1}a_j)^u(a_jz_{n+1})_s(q^{1-u+s}a_jz_{n+1})_u\\
&\quad=(-1)^u(q^{-s}/z_{n+1}a_j)^u(a_jz_{n+1})_{s+u}\end{split}\end{equation*}
if $u\in\{0,1\}$. Similar computations, eventually using \eqref{bc},
  reveal that \eqref{su} equals 
\begin{equation}\label{se}(-1)^u\frac{(q^{s+u})_t}{(q)_{s+t+u}}
\frac{\prod_{j=1}^{m+n+1}(a_jz_{n+1})_{s+u}}
{\prod_{j=1}^n(z_{n+1}/z_j)_{s+u}}
\prod_{j=1}^m\frac{(q^{s+u-1}w_jz_{n+1})_{y_j}}
{(w_jz_{n+1})_{y_j+s+u}}\end{equation}
if $u\in\{0,1\}$. Thus, the inner sum in \eqref{co} is proportional
to
$$\sum_{u=0}^{\min(1,N-t-|y|)}(-1)^u=\delta_{N-t-|y|,0}, $$
so we may assume $t+|y|=N$, $s=u=0$. But then the factor $(q^{s+u})_t$
in \eqref{se} equals zero unless $t=0$, so \eqref{co} is reduced to
the sum
$$\sum_{y_1+\dots+y_m=N}\frac{\De(\tilde wq^{(y,0,0)})}{\De(\tilde w)}
\prod_{k=1}^m\frac{\prod_{j=1}^{m+n+2}(w_k/\tilde a_j)_{y_k}}
{\prod_{j=1}^n(w_kz_j)_{y_k}\prod_{j=1}^{m+2}(qw_k/\tilde w_j)_{y_k}},
$$
where $s=0$ in \eqref{aw}.
It is easily verified that this simplifies to
$$\sum_{y_1+\dots+y_m=N}\frac{\De(wq^y)}{\De(w)}
\prod_{k=1}^m\frac{\prod_{j=1}^{m+n+1}(w_k/a_j)_{y_k}}
{\prod_{j=1}^{n+1}(w_kz_j)_{y_k}\prod_{j=1}^m(qw_k/w_j)_{y_k}}.$$
This completes the proof.
\end{proof}

\section{Further multiple Bailey transformations}

When $m=2$, the right-hand side of \eqref{te} has some additional symmetry
which allows us to obtain further transformations by iterating
Theorem \ref{t}. This was observed in \cite{ka2}, although the idea was not
fully exploited there. 

We first rewrite the case $m=2$ of \eqref{te} with $k=y_1=N-y_2$ as
summation index on the right-hand side, giving 
\begin{equation}\label{sm}\begin{split}&
\sum_{y_1+\dots+y_n=N}
\frac{\De(zq^y)}{\De(z)}\prod_{k=1}^n\frac{\prod_{j=1}^{n+2}(a_jz_k)_{y_k}}
{\prod_{j=1}^2(w_jz_k)_{y_k}\prod_{j=1}^n(qz_k/z_j)_{y_k}}\\
&=\frac{\prod_{j=1}^{n+2}(w_2/a_j)_N}{(w_2/w_1)_N(q)_N\prod_{j=1}^n(w_2z_j)_N}
\sum_{k=0}^N\Bigg(\frac{\theta(q^{2k-N}w_1/w_2)}{\theta(q^{-N}w_1/w_2)}
\\
&\quad\times
\frac{(q^{-N}w_1/w_2)_k(q^{-N})_k\prod_{j=1}^{n+2}(w_1/a_j)_k\prod_{j=1}^n
(q^{1-N}/w_2z_j)_k}{(q)_k(qw_1/w_2)_k\prod_{j=1}^{n+2}(q^{1-N}a_j/w_2)_k
\prod_{j=1}^n(w_1z_j)_k}\,q^k\Bigg),
\end{split}\end{equation}
where 
\begin{equation}\label{sc}
w_1w_2=z_1\dotsm z_na_1\dotsm a_{n+2}.\end{equation}

 Now let
$m\in\{0,1,\dots,n\}$, and write
\begin{equation}\label{xb}\begin{split}x&=(q^{1-N}a_1/w_1w_2,
    \dots,q^{1-N}a_{n-m}/w_1w_2,z_{n-m+1},\dots,z_n),\\
b&=(q^{N-1}w_1w_2z_{1},\dots,q^{N-1}w_1w_2z_{n-m},a_{n-m+1},\dots,a_{n+2}).
\end{split}\end{equation}
If we apply \eqref{sm} to the sum
$$\sum_{y_1+\dots+y_n=N}
\frac{\De(xq^y)}{\De(x)}\prod_{k=1}^n\frac{\prod_{j=1}^{n+2}(b_jx_k)_{y_k}}
{\prod_{j=1}^2(w_jx_k)_{y_k}\prod_{j=1}^n(qx_k/x_j)_{y_k}},$$
which is consistent with condition \eqref{sc}, we obtain the same
right-hand side up to the multiplier
\begin{equation*}\begin{split}
\frac{\prod_{j=1}^{n+2}(w_2/a_j)_N\prod_{j=1}^n(w_2x_j)_N}
{\prod_{j=1}^n(w_2z_j)_N\prod_{j=1}^{n+2}(w_2/b_j)_N}
&=\prod_{j=1}^{n-m}\frac{(w_2/a_j)_N(q^{1-N}a_j/w_1)_N}
{(w_2z_j)_N(q^{1-N}/w_1z_j)_N}\\
&=\prod_{j=1}^{n-m}\frac{(a_jz_j)^N(w_1/a_j)_N(w_2/a_j)_N}
{(w_1z_j)_N(w_2z_j)_N}.\end{split}\end{equation*}
This proves the following result.

\begin{corollary}\label{fc}
Assuming $w_1w_2=z_1\dotsm z_na_1\dotsm a_{n+2} $ and $0\leq m\leq n$,
we have
\begin{equation*}\begin{split}&
\sum_{\substack{y_1,\dots,y_n\geq 0\\y_1+\dots+y_n=N}}
\frac{\De(zq^y)}{\De(z)}\prod_{k=1}^n\frac{\prod_{j=1}^{n+2}(a_jz_k)_{y_k}}
{\prod_{j=1}^2(w_jz_k)_{y_k}\prod_{j=1}^n(qz_k/z_j)_{y_k}}\\
&=\prod_{j=1}^{n-m}\frac{(a_jz_j)^N(w_1/a_j)_N(w_2/a_j)_N}
{(w_1z_j)_N(w_2z_j)_N}\\
&\quad\times\sum_{\substack{y_1,\dots,y_n\geq 0\\y_1+\dots+y_n=N}}
\frac{\De(xq^y)}{\De(x)}\prod_{k=1}^n\frac{\prod_{j=1}^{n+2}(b_jx_k)_{y_k}}
{\prod_{j=1}^2(w_jx_k)_{y_k}\prod_{j=1}^n(qx_k/x_j)_{y_k}},
\end{split}\end{equation*}
where $x$ and $b$ are given by \eqref{xb}.
\end{corollary}

In the one-variable case, $n=2$, there are three choices of $m$: $m=2$,
which is trivial, $m=1$, which gives \eqref{fb} and $m=0$, which gives
\eqref{sb}. For general $n$ we have a sequence of non-trivial
identities: $m=0,1,\dots n-1$. 
The  case  $m=n-1$ is equivalent to
\cite[Corollorary 8.2]{r2}, which is an elliptic analogue of a
multiple Bailey transformation of Milne and Newcomb \cite{mn} (a
closely related identity was obtained in \cite{dg}).
The remaining identities, with $m\leq n-2$, appear to be new also in
the  non-elliptic case.  The  
extreme case $m=0$ is particularly elegant, so we
write it out explicitly. It gives a
 multivariable generalization of \eqref{sb} that is
 different from Theorem \ref{t}. 
 We have made the replacements $a_{n+1}=b$,
$a_{n+2}=c$, $w_1=d$, $w_2=e$.

\begin{corollary}\label{cs}
 Assuming $de=a_1\dotsm a_nbcz_1\dotsm z_n$, we have
\begin{equation*}\begin{split}&
\sum_{\substack{y_1,\dots,y_n\geq 0\\y_1+\dots+y_n=N}}
\frac{\De(zq^y)}{\De(z)}\prod_{j,k=1}^n\frac{(a_jz_k)_{y_k}}
{(qz_k/z_j)_{y_k}}\prod_{k=1}^n\frac{(bz_k)_{y_k}(cz_k)_{y_k}}
{(dz_k)_{y_k}(ez_k)_{y_k}}
\\
&=\prod_{j=1}^n\frac{(a_jz_j)^N(d/a_j)_N(e/a_j)_N}{(dz_j)_N(ez_j)_N}\\
&\quad\times\sum_{\substack{y_1,\dots,y_n\geq 0\\y_1+\dots+y_n=N}}
\frac{\De(aq^y)}{\De(a)}\prod_{j,k=1}^n\frac{(z_ja_k)_{y_k}}
{(qa_k/a_j)_{y_k}}\prod_{k=1}^n\frac{(q^{1-N}a_kb/de)_{y_k}
(q^{1-N}a_kc/de)_{y_k}}{(q^{1-N}a_k/d)_{y_k}(q^{1-N}a_k/e)_{y_k}}.
\end{split}\end{equation*}
\end{corollary}

Finally, we  give a companion identity to Corollary \ref{cs}, with
the sum supported on a hyper-rectangle rather than a simplex. There
are similar companions to the other cases of Corollary \ref{fc}, 
but these are  more complicated to write down.

We first replace $n$ by $n+1$ in  Corollary \ref{cs}, and assume that
$a_j=q^{-m_j}/z_j$, $1\leq j\leq n$, where $m_j$ are non-negative
integers such that $|m|\leq N$. Then all terms with $y_k>m_k$ for some
$k$ vanish. Eliminating $y_{n+1}$ from both summations we obtain
\begin{equation*}\begin{split}&
\sum_{y_1,\dots,y_n=0}^{m_1,\dots,m_n}
\frac{\De(zq^y)}{\De(z)}\prod_{k=1}^n
\frac{q^{N-|y|}\theta(z_kq^{y_k}/z_{n+1}q^{N-|y|})}{\theta(z_k/z_{n+1})}
\prod_{j,k=1}^n\frac{(q^{-m_j}z_k/z_j)_{y_k}}
{(qz_k/z_j)_{y_k}}\\
&\quad\times\prod_{k=1}^n\frac{(a_{n+1}z_k,bz_k,cz_k)_{y_k}}
{(qz_k/z_{n+1},dz_k,ez_k)_{y_k}}
\prod_{j=1}^n\frac{(q^{-m_j}z_{n+1}/z_j)_{N-|y|}}
{(qz_{n+1}/z_j)_{N-|y|}}\frac{(a_{n+1}z_{n+1},bz_{n+1},cz_{n+1})_{N-|y|}}
{(q,dz_{n+1},ez_{n+1})_{N-|y|}}\\
&=\prod_{k=1}^n\frac{(dz_kq^{m_k},ez_kq^{m_k})_N}{(dz_k,ez_k)_N}
\frac{(d/a_{n+1},e/a_{n+1})_N}{(dz_{n+1},ez_{n+1})_N}(a_{n+1}z_{n+1}q^{-|m|})^N
\\
&\quad\times\sum_{y_1,\dots,y_n=0}^{m_1,\dots,m_n}
\frac{\De(q^{y-m}/z)}{\De(q^{-m}/z)}\prod_{k=1}^n
\frac{q^{N-|y|}\theta(q^{y_k-m_k}/z_ka_{n+1}q^{N-|y|})}
{\theta(q^{-m_k}/z_ka_{n+1})}\\
&\quad\times\prod_{j,k=1}^n\frac{(q^{-m_k}z_j/z_k)_{y_k}}
{(q^{1-m_k+m_j}z_j/z_k)_{y_k}}\prod_{k=1}^n\frac{(q^{-m_k}z_{n+1}/z_k,
q^{1-N-m_k}b/dez_k,q^{1-N-m_k}c/dez_k)_{y_k}}
{(q^{1-m_k}/z_ka_{n+1},q^{1-N-m_k}/dz_k,q^{1-N-m_k}/ez_k)_{y_k}}\\
&\quad\times\prod_{j=1}^n\frac{(a_{n+1}z_j)_{N-|y|}}
{(q^{1+m_j}a_{n+1}z_j)_{N-|y|}}\frac{(a_{n+1}z_{n+1},q^{1-N}a_{n+1}b/de,
q^{1-N}a_{n+1}c/de)_{N-|y|}}
{(q,q^{1-N}a_{n+1}/d,q^{1-N}a_{n+1}/e)_{N-|y|}},
\end{split}\end{equation*}
where $deq^{|m|}=a_{n+1}bcz_{n+1}$.

We observe that the right-hand side will look nicer after the change
of variables $m_k\mapsto y_k-m_k$. We do that and also manipulate the
Pochhammer symbols so that $N$ never appears as a subscript.
After a tedious but straight-forward computation we arrive at
\begin{equation*}\begin{split}&
\sum_{y_1,\dots,y_n=0}^{m_1,\dots,m_n}
\frac{\De(zq^y)}{\De(z)}\prod_{k=1}^n
\frac{\theta(z_kq^{y_k+|y|-N}/z_{n+1})}{\theta(z_kq^{-N}/z_{n+1})}
\prod_{j,k=1}^n\frac{(q^{-m_j}z_k/z_j)_{y_k}}
{(qz_k/z_j)_{y_k}}\prod_{k=1}^n\frac{(a_{n+1}z_k,bz_k,cz_k)_{y_k}}
{(qz_k/z_{n+1},dz_k,ez_k)_{y_k}}\\
&\quad\times
\prod_{j=1}^n\frac{(q^{-N}z_j/z_{n+1})_{|y|}}
{(q^{1-N+m_j}z_j/z_{n+1})_{|y|}}\frac
{(q^{-N},q^{1-N}/dz_{n+1},q^{1-N}/ez_{n+1})_{|y|}}
{(q^{1-N}/a_{n+1}z_{n+1},q^{1-N}/bz_{n+1},q^{1-N}/cz_{n+1})_{|y|}}\,q^{|y|}\\
&=\prod_{k=1}^n\frac{(q^{1-N}z_k/z_{n+1},q^Ndez_k/b,q^Ndez_k/c)_{m_k}}
{(dz_k,ez_k,q^{N-|m|+1}a_{n+1}z_k)_{m_k}}\\
&\quad\times\frac{(d/a_{n+1},e/a_{n+1},q^{-N})_{|m|}}{(q^{N-|m|}a_{n+1}z_{n+1},
q^{N-|m|}bz_{n+1},q^{N-|m|}cz_{n+1})_{|m|}}\left(\frac{q^{2N}de}{bc}\right)
^{|m|}\\
&\quad\times\sum_{y_1,\dots,y_n=0}^{m_1,\dots,m_n}
\frac{\De(zq^y)}{\De(z)}\prod_{k=1}^n
\frac{\theta(z_ka_{n+1}q^{N-|m|+|y|+y_k})}
{\theta(z_ka_{n+1}q^{N-|m|})}\\
&\quad\times\prod_{j,k=1}^n\frac{(q^{-m_j}z_k/z_j)_{y_k}}
{(qz_k/z_j)_{y_k}}\prod_{k=1}^n\frac{(a_{n+1}z_k,q^Ndz_k,q^{N}ez_k)_{y_k}}
{(qz_k/z_{n+1},q^{N}dez_k/b,q^{N}dez_k/c)_{y_k}}\\
&\quad\times\prod_{j=1}^n\frac{(q^{N-|m|}z_ja_{n+1})_{|y|}}
{(q^{1+m_j+N-|m|}z_ja_{n+1})_{|y|}}\frac{(q^{N-|m|}a_{n+1}z_{n+1},q/cz_{n+1},
q/bz_{n+1})_{|y|}}
{(q^{1+N-|m|},q^{1-|m|}a_{n+1}/d,q^{1-|m|}a_{n+1}/e)_{|y|}}\,q^{|y|}.
\end{split}\end{equation*}
In this computation the following
identity, which is equivalent to \cite[Equation (3.8)]{r2},
is useful:
$$\frac{\De(1/z)}{\De(q^{-m}/z)}\prod_{j,k=1}^n\frac{(q^{-m_k}z_j/z_k)_{m_k}}
{(q^{1-m_k+m_j}z_j/z_k)_{m_k}}=(-1)^{|m|}q^{-|m|-\binom{|m|}{2}}. $$

To make  the connection with \eqref{sb}  transparent we
 make the change of parameters
$$(q^{-N},a_{n+1},b,c,d,e,z_{n+1})\mapsto(b,g,e,f,aq/c,aq/d,b/a). $$
We then obtain the following transformation, in the special case when
$b=q^{-N}$ with $N\geq |m|$ a non-negative integer.

\begin{corollary}\label{c3}
 Assuming $a^3q^{|m|+2}=bcdefg$, the following
  identity holds:
\begin{equation*}\begin{split}&
\sum_{y_1,\dots,y_n=0}^{m_1,\dots,m_n}
\frac{\De(zq^y)}{\De(z)}\prod_{k=1}^n
\frac{\theta(az_kq^{y_k+|y|})}{\theta(az_k)}
\prod_{j=1}^n\frac{(az_j)_{|y|}}
{(q^{1+m_j}az_j)_{|y|}}\frac
{(b,c,d)_{|y|}}{(aq/e,aq/f,aq/g)_{|y|}}\,q^{|y|}\\
&\quad\times
\prod_{j,k=1}^n\frac{(q^{-m_j}z_k/z_j)_{y_k}}
{(qz_k/z_j)_{y_k}}\prod_{k=1}^n\frac{(ez_k,fz_k,gz_k)_{y_k}}
{(aqz_k/b,aqz_k/c,aqz_k/d)_{y_k}}\\
&=g^{|m|}q^{-\sum_{j<k}m_jm_k}\frac{(b,aq/cg,aq/dg)_{|m|}}
{(aq/e,aq/f,aq/g)_{|m|}}\\
&\quad\times
\prod_{k=1}^n
\frac{z_k^{m_k}(aqz_k,q^{1+|m|-m_k}a/z_keg,q^{1+|m|-m_k}a/z_kfg)_{m_k}}
{(aqz_k/c,aqz_k/d,q^{|m|-m_k}b/gz_k)_{m_k}}\\
&\quad\times\sum_{y_1,\dots,y_n=0}^{m_1,\dots,m_n}
\frac{\De(zq^y)}{\De(z)}\prod_{k=1}^n
\frac{\theta(gz_kq^{y_k+|y|-|m|}/b)}{\theta(gz_kq^{-|m|}/b)}\\
&\quad\times
\prod_{j=1}^n\frac{(gz_jq^{-|m|}/b)_{|y|}}
{(gz_jq^{m_j+1-|m|}/b)_{|y|}}\frac
{(q^{-|m|}g/a,aq/be,aq/bf)_{|y|}}
{(q^{-|m|}cg/a,q^{-|m|}dg/a,q^{1-|m|}/b)_{|y|}}\,q^{|y|}\\
&\quad\times
\prod_{j,k=1}^n\frac{(q^{-m_j}z_k/z_j)_{y_k}}
{(qz_k/z_j)_{y_k}}\prod_{k=1}^n\frac{(aqz_k/bc,aqz_k/bd,gz_k)_{y_k}}
{(aqz_k/b,q^{-|m|}egz_k/a,q^{-|m|}fgz_k/a)_{y_k}}.
\end{split}\end{equation*}
\end{corollary}

To complete the proof we must extend the result from the case
$b=q^{-N}$ to generic $b$. This may be done exactly
 as  in the proof of Corollary 5.3
of \cite{r2}.
 That is, one considers the function $f(b)=L-R$, where $L$ and $R$ are
 the left- and right-hand side of the identity we want to prove, and
 where $c=a^3q^{|m|+2}/bdefg$ is viewed as depending on $b$ while the
 other parameters are fixed. It is then straight-forward to check that
 $f(pb)=f(b)$, where $p$ is the elliptic nome as in \eqref{th}.
Thus, $f(p^kq^{-N})=0$ for $k\in\mathbb Z$ and $N\in\mathbb Z_{\geq
  |m|}$. 
This is
  enough to conclude, by
 analytic continuation, that $f$ is identically zero.

\end{document}